\documentclass[11pt]{amsart}

\usepackage{amssymb}

\def\to{\mathchoice
{\longrightarrow}
{\rightarrow}
{\rightarrow}
{\rightarrow}}

\def\mapsto{\mathchoice
{\DOTSB\mapstochar\longrightarrow}
{\DOTSB\mapstochar\rightarrow}
{\DOTSB\mapstochar\rightarrow}
{\DOTSB\mapstochar\rightarrow}}

\def\hookrightarrow{\mathchoice
{\DOTSB\lhook\joinrel\relbar\joinrel\rightarrow}
{\DOTSB\lhook\joinrel\rightarrow}
{\DOTSB\lhook\joinrel\rightarrow}
{\DOTSB\lhook\joinrel\rightarrow}}

\swapnumbers

\newtheorem{theorem}{Theorem}[section]
\newtheorem*{theorem*}{Theorem}
\newtheorem{proposition}[theorem]{Proposition}
\newtheorem{lemma}[theorem]{Lemma}
\newtheorem{corollary}[theorem]{Corollary}

\theoremstyle{definition}

\newtheorem{definition}[theorem]{Definition}

\newcommand{\C}{\mathbb{C}}
\newcommand{\Q}{\mathbb{Q}}
\newcommand{\F}{\mathbb{F}}

\newcommand{\eps}{\varepsilon}
\renewcommand{\o}{\otimes}

\DeclareMathOperator{\Ind}{Ind}

\newcommand{\p}{\partial}

\newcommand{\<}{\langle}
\renewcommand{\>}{\rangle}

\newcommand{\bull}{\bullet}
\renewcommand{\SS}{\mathbb{S}}

\DeclareMathOperator{\GL}{GL}
\DeclareMathOperator{\gl}{\mathfrak{gl}}
\DeclareMathOperator{\SP}{Sp}
\newcommand{\half}{\tfrac12}
\newcommand{\HH}{\mathsf{H}}
\newcommand{\Lie}{\mathsf{Lie}}
\newcommand{\A}{\mathsf{A}}
\newcommand{\B}{\mathsf{B}}
\renewcommand{\L}{\mathsf{L}}
\newcommand{\K}{\mathsf{K}}
\newcommand{\DD}{\mathcal{D}}
\DeclareMathOperator{\Wedge}{\Lambda}
\DeclareMathOperator{\Hom}{Hom}
\DeclareMathOperator{\Schur}{\mathsf{S}}
\newcommand{\CO}{\mathcal{O}}
\newcommand{\CP}{\mathcal{P}}
\newcommand{\FF}{\mathsf{F}}
\newcommand{\Mbar}{\overline{\mathcal{M}}}
\newcommand{\CM}{\mathcal{M}}
\renewcommand{\.}{\wedge}

\begin{document}

\title{The homology groups of some two-step nilpotent Lie algebras
associated to symplectic vector spaces}

\author{E. Getzler}

\address{Department of Mathematics, Northwestern University, Evanston, IL
60208-2730, USA}

\email{getzler@math.nwu.edu}

\maketitle

Let $\Phi:V\mapsto\Phi(V)$ be a polynomial functor from the category of
vector spaces (over a field $\F$ of characteristic zero) to the category of
Lie algebras. In this paper, we study the functors $H_k(\Phi):V\mapsto
H_k(\Phi(V))$ from vector spaces to vector spaces obtained by composing
$\Phi$ with the $k$th Lie algebra homology group functor. These functors
are also polynomial functors, and are best studied by expressing them as
explicit Schur functors.

The simplest example is obtained by taking $\Phi$ to be the identity
functor, which assigns to a vector space $V$ the Lie algebra $V$ with
vanishing bracket. In this case, $H_k(\Phi)$ is the $k$th exterior power.

A more complicated example was investigated by Sigg \cite{Sigg}. Take
$\Phi$ to be the free $2$-step nilpotent Lie algebra functor
$V\mapsto\Lie_2(V)=V\oplus\Wedge^2V$, with bracket
$$
[(v_1,a_1),(v_2,a_2)] = (0,v_1\wedge v_2) , \quad\quad v_i\in V,
a_i\in\Wedge^2V .
$$
If $\lambda$ is a Young diagram, let $\Schur^\lambda$ be the Schur functor
associated to $\lambda$ (cf.\ Fulton-Harris \cite{FH}); in particular,
$\Schur^{(1^k)}$ is the $k$th exterior power, while $\Schur^{(k)}$ is the
$k$th symmetric power. Let $\lambda^*$ be the conjugate partition of
$\lambda$, defined by $\lambda_i^*=\sup\{j\mid \lambda_j\ge i\}$. Introduce
the set $\CO_k$ of Young diagrams such that $\lambda$ is self-conjugate,
$\lambda=\lambda^*$, and $2k = |\lambda| + \sup \{ i \mid \lambda_i \ge
i\}$. Sigg proves that
$$
H_k(\Lie_2) \cong \sum_{\lambda\in\CO_k} \Schur^\lambda .
$$
For example, $H_1(\Lie_2)$ is the identity functor,
$H_2(\Lie_2)\cong\Schur^{(2,1)}$, and
$H_3(\Lie_2)\cong\Schur^{(3,1^2)}\oplus\Schur^{(2^2)}$.

In this paper, we prove an analogue of Sigg's result. Let $\HH$ be the
symplectic vector space $\F\oplus\F$, with symplectic form
$\<(a,b),(c,d)\>=ad-bc$. Let $\L_\HH(V)$ be the Lie algebra $(\HH\o
V)\oplus \Schur^2(V)$, with bracket
$$
[(v_1,w_1;a_1),(v_2,w_2;a_2)] = (0,0;v_1\cdot w_2-v_2\cdot w_1) ,
\quad\quad (v_i,w_i)\in\HH\o V, a_i\in \Schur^2(V) .
$$
The homology $H_k(\L_\HH)$ is more complicated than that of $\Lie_2$, even
when $\HH$ is two-dimensional, and we have not been able to calculate it
completely. As an illustration,
$$
H_k(\L_\HH) \cong \begin{cases} \bigl( \HH_0 \o \Schur^{(0)} \bigr), & k=0
, \\
\bigl( \HH_1 \o \Schur^{(1)} \bigr) , & k=1 , \\
\bigl( \HH_2 \o \Schur^{(1^2)} \bigr) \oplus \bigl( \HH_1 \o \Schur^{(3)}
\bigr) , & k=2 ,\\
\bigl( \HH_3 \o \Schur^{(1^3)} \bigr) \oplus \bigl( \HH_2 \o \Schur^{(3,1)}
\bigr) \oplus \bigl( \HH_0 \o \Schur^{(4)} \bigr) & k=3 ,
\end{cases}$$
where $\HH_k$ is the $k$th symmetric power of $\HH$ (and is thus
$k+1$-dimensional).

The Lie algebras $\L_\HH(V)$ satisfy Poincar\'e duality, since their
associated simply connected Lie group is contractible and contains a
cocompact lattice. (We owe this remark to P. Etingof.) For example,
$\Schur^\lambda(\F)$ is nonzero only if $\lambda$ has length $1$; we see
that
$$
H_\bull(\L_\HH(\F)) \cong \F \oplus \HH[1] \oplus \HH[3] \oplus \F[4] ,
$$
where $V[k]$ is the vector space $V$ shifted into degree $k$. In higher
dimensions, Poincar\'e duality is difficult to see directly.

As we explain in \cite{genus1}, the homology groups $H_k(\L_\HH)$ are
closely related to the $E_2$-terms of the Leray-Serre spectral sequence for
the fibrations $\CM_{g,n}\to\CM_g$ (or, in genus $1$, of the fibrations
$\CM_{1,n}\to\CM_{1,1}$). In particular, the summand $\HH_0 \o
\Schur^{(4)}$ of $H_3(\L_\HH)$ gives rise to the relation in
$H^4(\Mbar_{1,4},\Q)$ discovered in \cite{elliptic}.

\subsection*{Acknowledgements} We are grateful to M. Kapranov for informing
us of Sigg's work. This research is partially supported by the NSF.

\section{Lie $\SS$-algebras}

In this section, we recall parts of the formalism of operads, referring to
Getzler-Jones \cite{GJ} for further details. This formalism is closely
related to Joyal's theory of species and analytic functors (Joyal
\cite{Joyal}).

\begin{definition}
An $\SS$-module is a functor from $\SS$, the groupoid formed by taking the
union of the symmetric groups $\SS_n$, $n\ge0$, to the category of vector
spaces.
\end{definition}

Associated to an $\SS$-module $\A$ is the functor from the category of
vector spaces to itself,
$$
V \mapsto \A(V) = \sum_{k=0}^\infty \bigl( \A(k) \o V^{\o k}
\bigr)_{\SS_k} .
$$
This is a generalization of the notion of a Schur functor, which is the
special case where $\A$ is an irreducible representation of $\SS_n$.

\begin{definition}
A polynomial functor $\Phi$ is a functor from the category of vector spaces
to itself such that the map $\Phi:\Hom(V,W)\to\Hom(\Phi(V),\Phi(W))$ is
polynomial for all vector spaces $V$ and $W$. An analytic functor $\Phi$ is
a direct image of polynomial maps.
\end{definition}

To an analytic functor $\Phi$, we may associate the $\SS$-module
$$
\A(n) = \Phi(\F x_1\oplus\dots\oplus\F x_n)_{(1,\dots,1)} \subset \Phi(\F
x_1\oplus\dots\oplus\F x_n) ,
$$
the summand of $\Phi(\F x_1\oplus\dots\oplus\F x_n)$ homogeneous of degree
$1$ in each of the generators $x_i$. We call $\A$ the $\SS$-module of
Taylor coefficients of $\Phi$. The following theorem is proved in
Appendix~A of Macdonald \cite{Macdonald}.
\begin{theorem}
There is an equivalence of categories between the category of $\SS$-modules
and the category of analytic functors: to an $\SS$-module, we associate the
functor $V\mapsto\A(V)$, while to an analytic functor $\Phi$, we
associate its $\SS$-module of Taylor coefficients.
\end{theorem}

Any $\SS$-module $\A$ extends to a functor on the category of finite sets
and bijections: if $S$ is a finite set of cardinality $n$, we have
$$
\A(S) = \biggl( \sum_{\substack{f:[n]\to S \\ \text{bijective}}} \A(n)
\biggr)_{\SS_n} ,
$$
where $[n]=\{1,\dots,n\}$. The category of $\SS$-modules has a monoidal
structure, defined by the formula
$$
(\A\circ\B)(n) = \sum_{k=0}^\infty \biggl( \A(k) \o \sum_{f:[n]\to[k]}
\B(f^{-1}(1)) \o \dots \o \B(f^{-1}(k)) \biggr)_{\SS_k} .
$$
This definition is motivated by the composition formula
$(\A\circ\B)(V)\cong\A(\B(V))$.

\begin{definition}
An \textbf{operad} is a monoid in the category of $\SS$-modules, with
respect to the above monoidal structure.
\end{definition}

We see that the structure of an operad on an $\SS$-module $\A$ is the same
as the structure of a triple on the associated analytic functor
$V\mapsto\A(V)$.

The Lie operad $\Lie$ is the operad whose associated analytic functor is
the functor taking a vector space to its free Lie algebra.
\begin{definition}
A Lie $\SS$-algebra $\L$ is a left $\Lie$-module in the category of
$\SS$-modules.
\end{definition}

Lie $\SS$-algebras are essentially the same things as analytic functors
from the category of vector spaces to the category of Lie algebras; more
precisely, they are the collections of Taylor coefficients of such
functors.

If we unravel the definition of a Lie $\SS$-algebra, we see that it is an
$\SS$-module $\L$ with $\SS_n$-equivariant brackets
$$
[-,-] : \Ind^{\SS_n}_{\SS_k\times\SS_{n-k}} \L(k) \o \L(n-k) \to \L(n)
$$
for $0\le k\le n$, such that if $a_i\in\L(n_i)$, $i=1,2,3$, the following
expressions vanish:
\begin{gather*}
[a_1,a_2] - [a_2,a_1] \in \Ind^{\SS_n}_{\SS_{n_1}\times\SS_{n_2}} \bigl(
\L(n_1) \o \L(n_2) \bigr) \quad \text{and} \\
[a_1,[a_2,a_3]] + [a_2,[a_3,a_1]] + [a_3,[a_1,a_2]] \in
\Ind^{\SS_n}_{\SS_{n_1}\times\SS_{n_2}\times\SS_{n_3}} \bigl( \L(n_1) \o
\L(n_2) \o \L(n_3) \bigr) .
\end{gather*}

If $L$ is a Lie algebra, let $\K_\bull(L)$ be the Chevalley-Eilenberg
complex of $L$. Recall that $\K_k(L)=\Wedge^kL$ is the $k$th exterior power
of $L$, and the differential $\p:\K_k(L)\to\K_{k-1}(L)$ is given by the
formula
$$
\p (a_1\.\dots\.a_k) = \sum_{1\le i<j\le k} (-1)^{i-j+1} [a_i,a_j] \. a_1
\. \dots \. \widehat{a_i} \. \dots \widehat{a_j} \. \dots \. a_k .
$$
If $\L$ is a Lie $\SS$-algebra, we obtain a sequence of analytic functors
$V\mapsto\bigl( \K_\bull(\L(V)) , \p \bigr)$. Define the
Chevalley-Eilenberg complex of the Lie $\SS$-algebra $\L$ to be the Taylor
coefficients of this complex of analytic functors. In other words,
$\K_k(\L)=\Wedge^k\circ\L$, where $\Wedge^k$ is the $\SS$-module
$$
\Wedge^k(n) = \begin{cases} \Schur^{(1^k)} , & k=n , \\ 0 , & k\ne n
. \end{cases}
$$
The differential $\p:\K_k(\L(V))\to\K_{k-1}(\L(V))$ is a natural
transformation of analytic functors, and hence induces a map of
$\SS$-modules $\p:\K_k(\L)\to\K_{k-1}(\L)$. Clearly, we have $\p^2=0$.
\begin{definition}
The $k$th homology group $H_k(\L)$ of the Lie $\SS$-algebra $L$ is the
$k$th homology group of the complex of $\SS$-modules $(\K_\bull(\L),\p)$.
\end{definition}

Thus, $H_k(\L)$ is an $\SS$-module for each $k\ge0$.

\section{Examples of Lie $\SS$-algebras}

As a left module over itself, the Lie operad $\Lie$ is a Lie $\SS$-algebra;
the corresponding analytic functor is the free Lie algebra functor. More
generally, define $\Lie_d$, $1\le d\le\infty$, by
$$
\Lie_d(n) = \begin{cases} \Lie(n) , & n\le d , \\ 0 , & n>d . \end{cases}
$$
Each of these is a Lie $\SS$-algebra; the brackets
$\Lie_d(k)\o\Lie_d(n-k)\to\Lie_d(n)$ are defined as for $\Lie$ if $n\le d$,
and of course vanish if $n>d$. The analytic functor associated to the Lie
$\SS$-module $\Lie_d$ is known as the free $d$-step nilpotent Lie algebra.
We may view Sigg's theorem \cite{Sigg} as the calculation of the homology
of the Lie $\SS$-algebra $\Lie_2$:
$$
H_k(\Lie_2)(n) \cong \sum_{\{\lambda\in\CO_k\mid|\lambda|=n\}}
\Schur^\lambda .
$$
Here, we use the same notation for the representation of the symmetric
group $\SS_n$ with the Young diagram $\lambda$ as for the associated Schur
functor $\Schur^\lambda$.

The tensor product $R\o\L$ of a Lie $\SS$-algebra $\L$ with a commutative
algebra $R$ is again a Lie $\SS$-algebra. For example, let $M$ be a
differentiable manifold and let $\Omega^\bull(M)$ be the differential
graded algebra of complex differential forms. The homology of differential
graded Lie $\SS$-algebras is defined in a manner analogous to the
definition of the homology of Lie $\SS$-algebras, except that we must add
to the Chevalley-Eilenberg differential $\p$ the internal differential $d$
in defining the homology groups. Let $\FF(M,n)$ be the $n$th configuration
space of $M$, defined by
$$
\FF(M,n) = \{ i : [n] \to  M \mid \text{$i$ is an embedding} \} .
$$
Let $j(n):\FF(M,n)\to M^n$ be the open embedding of the configuration
space.  The resolution of the sheaf $j(n)_!j(n)^*\C$ on $M^n$ constructed
in \cite{config2} may be identified with the twist of the
Chevalley-Eilenberg complex $\K_\bull(\Omega^\bull(M)\o\Lie)(n)$ by the
alternating character $\eps(n)$ of $\SS_n$. This yields natural
isomorphisms
$$
H^\bull(\FF(M,n),\C)[n] \cong H_\bull(\Omega^\bull(M)\o\Lie)(n) \o \eps(n)
.
$$
In particular, if $M$ is a compact manifold whose cohomology over $\C$ is
formal (such as a compact K\"ahler manifold), we see that
$$
H^\bull(\FF(M,n),\C)[n] \cong H_\bull(H^\bull(M,\C)\o\Lie)(n) \o \eps(n) .
$$
This reformulates a theorem of Totaro \cite{Totaro}.

Another example of a Lie $\SS$-algebra is associated to a symplectic vector
space $\HH$ with symplectic form $\<-,-\>$: set $\L_\HH(1)=\HH$, and let
$\L_\HH(2)$ be the trivial representation $\Schur^{(2)}$ of $\SS_2$. The
Chevalley-Eilenberg complex of $\L_\HH$ is familiar from Weyl's
construction of the irreducible representations of the symplectic group
$\SP(\HH)$: we have
$$
\K_n(\L_\HH)(n+\ell) = \begin{cases}
\Ind_{\SS_\ell\wr\SS_2\times\SS_{n-\ell}}^{\SS_{n+\ell}} \Bigl( \bigl(
\Schur^{(2)} \bigr)^{\o\ell} \o \Schur^{(1^{n-\ell})} \Bigr) \o
\HH^{\o(n-\ell)} , & \ell\ge0 , \\ 0 , & \ell<0 .
\end{cases}
$$
In particular, $\K_n(\L_\HH)(n)\cong\Schur^{(1^n)}\o\HH^{\o n}$, and
$$
\K_n(\L_\HH)(n+1) \cong \sum_{1\le i<j\le n} \Schur^{(1^{n-1})} \o
\HH^{\o(n-1)} \o x_{ij} .
$$
The differential $\p:\K_n(\L_\HH)(n)\to\K_{n+1}(\L_\HH)(n)$ is given by
$$
\p ( e_1 \o \dots \o e_n ) = \sum_{1\le i<j\le n} (-1)^{j-i+1} \<e_i,e_j\>
\, e_1 \o \dots \o \widehat{e_i} \o \dots \o \widehat{e_j} \o \dots \o e_n
\o x_{ij} .
$$
If $\Schur^{\<\lambda\>}(\HH)$ is the irreducible representation of
$\SP(\HH)$ associated to the Young diagram $\lambda$, it follows that
$$
H_n(\L_\HH)(n) \cong \sum_{|\lambda|=n} \Schur^{\<\lambda\>}(\HH) \o
\Schur^{\lambda^*} .
$$
For example, if $\dim(\HH)=2$, denoting the $k$th symmetric power
$\Schur^{\<k\>}(\HH)$ of $\HH$ by $\HH_k$, we have
$$
\K_n(\L_\HH)(n) \cong \sum_{j=0}^{[\frac{n}{2}]} \HH_{n-2j} \o
\Schur^{(2^j,1^{n-2j})} ,
$$
and $H_n(\L_\HH)(n) \cong \HH_n \o \Schur^{(1^n)}$.

\section{The Chevalley-Eilenberg complex of $\L_\HH$}

We now turn to the closer study of the Chevalley-Eilenberg complex of the
Lie $\SS$-algebra $\L_\HH$. To this end, choose a basis $\{e_a \mid 1\le
a\le 2g \}$ for $\HH$, with symplectic form
$$
\<e_a,e_b\>=\eta_{ab} .
$$
Let $\eta^{ab}$ be the inverse matrix to $\eta_{ab}$:
$$
\sum_{b=1}^{2g} \eta^{ab} \eta_{bc} = \delta^a_c .
$$

Let $V$ be a vector space with basis $\{E_i\mid 1\le i\le r\}$; the
symmetric square $\Schur^2(V)$ has basis $\{E_{ij}=E_iE_j\mid 1\le i\le j
\le r\}$.

The nilpotent Lie algebra $\L_\HH(V)=(\HH\o V)\oplus\Schur^2(V)$ has centre
$\Schur^2(V)$, and the restriction of its Lie bracket to $\HH\o V$ is
$$
[ e_a\o E_i , e_b\o E_j ] = \eta_{ab} E_{ij} .
$$
The Chevalley-Eilenberg complex of $\L_\HH(V)$ is the graded vector space
$\Wedge^\bull(\HH\o V) \o \Wedge^\bull(\Schur^2(V))$. Denote by $\eps_i^a$
the operation of exterior multiplication by $e_a\o E_i$ on this complex,
and let $\iota_a^i$ be its adjoint, characterized by the (graded)
commutation relations
$$
[\iota_a^i,\eps_j^b] = \delta_j^i \delta^b_a .
$$
Let $\eps_{ij}=\eps_{ji}$ be the operation of exterior multiplication by
$E_{ij}$ on the Chevalley-Eilenberg complex, and let $\iota^{ij}$ be its
adjoint, characterized by the commutation relations
\begin{equation} \label{commute.ij}
[\iota^{ij},\eps_{kl}] = \delta_k^i \delta_l^j + \delta_l^i \delta_k^j .
\end{equation}

The differential $\p$ of the Chevalley-Eilenberg complex and its adjoint
$\p^*$ are given by the formulas
$$
\p = \half \sum_{i,j,a,b} \eta^{ab} \eps_{ij} \iota^i_a \iota^j_b , \quad
\p^* = - \half \sum_{i,j,a,b} \eta_{ab} \eps^a_i \eps^b_j \iota^{ij} .
$$
The following theorem is the most powerful idea in the calculation of the
cohomology of nilpotent Lie algebras.
\begin{theorem}[Kostant \cite{Kostant}]
The kernel of the Laplacian $\Delta=[\p^*,\p]$ on the Chevalley-Eilenberg
complex is isomorphic to the homology of the Lie algebra $\L_\HH(V)$.
\end{theorem}

Sigg \cite{Sigg} has calculated the Laplacian $\Delta$ for the free
$2$-step nilpotent Lie algebra $\Lie_2(V)=V\oplus\Wedge^2V$. Our
calculation is modelled on his, with some modifications brought on by the
introduction of the symplectic vector space $\HH$.

The complexity of our notation is reduced by adopting the Einstein
summation convention: indices $i,j,\dots$ lie in the set $\{1,\dots,r\}$,
indices $a,b,\dots$ in the set $\{1,\dots,2g\}$, and we sum over repeated
pairs of indices if one is a subscript and one is a superscript.
\begin{lemma} \label{Laplacian1}
$\Delta = \eps_{ij}\eps_k^a\iota^i_a\iota^{jk} - \half \eta_{ab} \eta^{cd}
\eps_i^a \eps_j^b \iota^i_c \iota^j_d - g \, \eps_{ij}\iota^{ij}$
\end{lemma}
\begin{proof}
We have
$$
4 [\p^*,\p] = - [ \eta_{ab} \eps^a_i \eps^b_j \iota^{ij} , \eta^{cd}
\eps_{kl} \iota^k_c \iota^l_d] = - \eta_{ab} \eta^{cd} \eps_{kl} [ \eps^a_i
\eps^b_j , \iota^k_c \iota^l_d ] \iota^{ij} - \eta_{ab} \eta^{cd} \eps^a_i
\eps^b_j [ \iota^{ij} , \eps_{kl} ] \iota^k_c \iota^l_d .
$$
The first term of the right-hand side is calculated as follows,
\begin{align*}
- \eta_{ab} \eta^{cd} [ \eps^a_i \eps^b_j , \iota^k_c \iota^l_d ] &=
- \eta_{ab} \eta^{cd} \eps^a_i [ \eps^b_j , \iota^k_c \iota^l_d ]
- \eta_{ab} \eta^{cd} [ \eps^a_i , \iota^k_c \iota^l_d ] \eps^b_j \\
&=  \delta^k_j \eps^a_i \iota^l_a + \delta^l_j \eps^a_i \iota^k_a
- \delta^k_i \iota^l_a \eps^a_j - \delta^l_i \iota^k_a \eps^a_j \\
&= \delta^k_j \eps^a_i \iota^l_a + \delta^l_j \eps^a_i \iota^k_a
+ \delta^k_i \eps^a_j \iota^l_a + \delta^l_i \eps^a_j \iota^k_a - 2g \,
\delta^k_i \delta^l_j - 2g \, \delta^l_i \delta^k_j ,
\end{align*}
while the second term is calculated by \eqref{commute.ij}.
\end{proof}

\section{The Casimir operator of $\GL(V)$}

If $V$ is a vector space with basis $\{E_i\mid 1\le i\le n\}$, the Lie
algebra of $\GL(V)$ has basis $\{E_i^j\mid 1\le i,j\le n\}$, with
commutation relations
$$
[E_i^j,E_k^l] = \delta^j_k E_i^l - \delta_i^l E^j_k .
$$
The centre of $\GL(V)$ is spanned by $\DD = E_i^i$, and the Casimir
operator is the element of the centre of $U(\gl(V))$ given by the formula
$$
\Delta_{\GL(V)} = E_i^jE_j^i .
$$
Let $c_\lambda$ be the eigenvalue of the Casimir operator $\Delta_{\GL(V)}$
on the representation $\Schur^\lambda(V)$ of $\GL(V)$ with highest weight
vector $\lambda=(\lambda_1,\dots,\lambda_r)$. Since the sum of the positive
roots of $\GL(V)$ equals $2\rho = (2r-1,2r-3,\dots,3-2r,1-2r)$, the theory
of semisimple Lie algebras shows that, up to an overall factor,
\begin{equation} \label{casimir}
c_\lambda = \|\lambda\|^2 + 2(\rho,\lambda) = \sum_{i=1}^r
\lambda_i(\lambda_i+r-2i+1) .
\end{equation}
To see that this factor equals $1$, observe that on the fundamental
representation $V$, with highest weight $(1,0,\dots,0)$, the Casimir has
eigenvalue $r$.

Given a Young diagram $\lambda$, let
$$
n(\lambda) = \sum_{i\ge1} (i-1)\lambda_i = \sum_{i\ge1}
\binom{\lambda_i^*}{2} .
$$
\begin{lemma} \label{Casimir}
$c_\lambda = r|\lambda| + 2n(\lambda^*) - 2n(\lambda)
= \sum_{i=1}^\infty \lambda^*_i(r-\lambda^*_i+2i-1)$
\end{lemma}
\begin{proof}
The proof follows from rearranging \eqref{casimir}:
$$
c_\lambda = r|\lambda| + 2 \sum_{i=1}^r \binom{\lambda_i}{2} - 2
\sum_{i=1}^r (i-1) \lambda_i .
\qed$$\def\qed{}
\end{proof}

Recall the dominance order on Young diagrams:
$$
\text{$\lambda\ge\mu$ if $|\lambda|=|\mu|$ and
$\lambda_1+\dots+\lambda_i\ge\mu_1+\dots+\mu_i$ for all $i\ge1$.}
$$
If $\lambda\ge\mu$, then $\mu^*\ge\lambda^*$ (Macdonald, I.1.11
\cite{Macdonald}).
\begin{corollary} \label{dominant}
If $\lambda\ge\mu$, then $c_\lambda \ge c_\mu$, with equality only if
$\lambda=\mu$.
\end{corollary}
\begin{proof}
If $\lambda\ge\mu$, we have
$$
n(\lambda) = \sum_{i\ge1} (i-1)\lambda_i = \sum_{i\ge1} \sum_{j>i}
\lambda_i = \sum_{i\ge1} \Bigl( |\lambda| - \sum_{j=1}^i \lambda_i \Bigr)
\le \sum_{i\ge1} \Bigl( |\mu| - \sum_{j=1}^i \mu_i \Bigr) = n(\mu) .
$$
Likewise, $n(\lambda^*)\ge n(\mu^*)$. In both cases, equality holds only if
$\lambda=\mu$. The corollary now follows from Lemma~\ref{Casimir}.
\end{proof}

\begin{corollary} \label{Dominant}
On the tensor product $\Schur^\lambda(V)\o\Schur^\mu(V)$, the Casimir
operator $\Delta_{\GL(V)}$ is bounded above by
$c_\lambda+c_\mu+2(\lambda,\mu)$, with equality only on
$\Schur^{\lambda+\mu}(V)\hookrightarrow\Schur^\lambda(V)\o\Schur^\mu(V)$.
\end{corollary}
\begin{proof}
There can only be a nonzero morphism
$\Schur^\nu(V)\hookrightarrow\Schur^\lambda(V)\o\Schur^\mu(V)$ if
$\nu\le\lambda+\mu$. It follows from Corollary \ref{dominant} that
$$
c_\nu \le c_{\lambda+\mu} = \|\lambda+\mu\|^2 + 2 (\rho,\lambda+\mu) =
\|\lambda\|^2 + 2 (\rho,\lambda) + \|\mu\|^2 + 2 (\rho,\mu) +
2(\lambda,\mu) .
\qed$$\def\qed{}
\end{proof}

\section{A formula for the Laplacian}

In this section, we prove the following explicit formula for the Laplacian
$\Delta$ on the Chevalley-Eilenberg complex $\K_\bull(\L_\HH(V))$.
\begin{theorem} \label{main}
$\Delta = \half \bigl( \Delta_{\SP(\HH)} + \Delta_{\GL(V)} - (r+2g+1) \DD
\bigr)$
\end{theorem}

Theorem \ref{main} will follow by combining the results of
Lemmas~\ref{Laplacian1}, \ref{Laplacian2} and \ref{Laplacian3}. The Lie
algebra of $\GL(V)$ acts on $\K_\bull(\L_\HH(V))$ via the operations
$$
E_i^j = \eps_i^a \iota^j_a + \eps_{ik} \iota^{jk} .
$$
It follows that $\DD = \eps_i^a \iota^i_a + \eps_{ij} \iota^{ij}$, while
the Casimir operator for $\GL(V)$ acts on $\K_\bull(\L_\HH(V))$ as follows.
\begin{lemma} \label{Laplacian2}
$\Delta_{\GL(V)} = \eps_i^a\eps_j^b \iota^i_b\iota^j_a + 2 \,
\eps_{ij}\eps_k^a\iota^i_a\iota^{jk} + r \eps_i^a \iota^i_a + (r+1)
\eps_{ij} \iota^{ij}$
\end{lemma}
\begin{proof}
We have
\begin{align*}
E_i^jE_j^i & = ( \eps_i^a \iota^j_a + \eps_{ik} \iota^{jk} ) ( \eps_j^b
\iota^i_b + \eps_{jl} \iota^{il} ) \\ &= \eps_i^a \iota^j_a \eps_j^b
\iota^i_b + \eps_i^a \iota^j_a \eps_{jl} \iota^{il} + \eps_{ik} \iota^{jk}
\eps_j^b \iota^i_b + \eps_{ik} \iota^{jk} \eps_{jl} \iota^{il} \\ &= -
\eps_i^a \eps_j^b \iota^j_a \iota^i_b + r \eps_i^a \iota^i_a + \eps_{jl}
\eps_i^a \iota^j_a \iota^{il} + \eps_{ik} \eps_j^b \iota^i_b \iota^{jk} -
\eps_{ik} \eps_{jl} \iota^{jk} \iota^{il} + (r+1) \eps_{ij} \iota^{ij} .
\end{align*}
The (a)symmetries of $\eps_{ik} \eps_{jl} \iota^{jk} \iota^{il}$ force it
to vanish, and the result follows.
\end{proof}

The Lie algebra of $\GL(\HH)$ acts on the Chevalley-Eilenberg complex of
$\L_\HH(V)$ by the operators
$$
\{ e^a_b = \eps_i^a \iota^i_b \mid 1\le a,b\le 2g \} ,
$$
and the Lie subalgebra $\SP(\HH)\subset\GL(\HH)$ is spanned by the
operators
$$
\{ e_{ab}+e_{ba} \mid 1\le a\le b\le 2g \} ,
$$
where $e_{ab}=\eta_{ac}e^c_b$. The Casimir operator of $\SP(\HH)$ is given
by the formula
$$
\Delta_{\SP(\HH)} =- \half \eta^{ac} \eta^{bd} \bigl( e_{ab} + e_{ba}
\bigr) \bigl( e_{cd} + e_{dc} \bigr) = - \eta^{ac} \eta^{bd} e_{ab} e_{cd}
- \eta^{ac} \eta^{bd} e_{ab} e_{dc} .
$$
\begin{lemma} \label{Laplacian3}
$\Delta_{\SP(\HH)} = -\eps_i^a\eps_j^b\iota^i_b\iota^j_a
-\eta_{ab}\eta^{cd}\eps_i^a\eps_j^b\iota^i_c\iota^j_d +
(2g+1)\eps_i^a\iota^i_a$
\end{lemma}
\begin{proof}
We have
\begin{align*}
\eta^{ac} \eta^{bd} e_{ab} e_{cd} &= \eta^{ac} \eta^{bd} \eta_{aa'}
\eta_{cc'} \eps_i^{a'} \iota^i_b \eps_j^{c'} \iota^j_d = - \eta_{ac}
\eta^{bd} \eps_i^a \iota^i_b \eps_j^c \iota^j_d = \eta_{ac} \eta^{bd}
\eps_i^a \eps_j^c \iota^i_b \iota^j_d - \eps_i^a \iota^i_a \\
\eta^{ac} \eta^{bd} e_{ab} e_{dc} &= \eta^{ac} \eta^{bd} \eta_{aa'}
\eta_{dd'} \eps_i^{a'} \iota^i_b \eps_j^{d'} \iota^j_c = - \eps_i^a
\iota^i_b \eps_j^b \iota^j_a = \eps_i^a \eps_j^b \iota^i_b \iota^j_a - 2g
\eps_i^a \iota^i_a .
\qed\end{align*}
\def\qed{}
\end{proof}

\section{The case $g=1$}

In this section, we apply our results in the special case $g=1$, in which
the symplectic vector space $\HH$ is two-dimensional. Recall Frobenius's
notation for partitions: if $\alpha_1>\dots>\alpha_d\ge0$ and
$\beta_1>\dots>\beta_d\ge0$,
$$
(\alpha_1,\dots,\alpha_d|\beta_1,\dots,\beta_d)
$$
is the partition of $\alpha_1+\dots+\alpha_d+\beta_1+\dots+\beta_d+d$ whose
$i$th part equals $\alpha_i+i$ for $i\le d$, and $\sup\{j \mid \beta_j+j
\ge i\}$ for $i>d$. For example, $(\alpha|\beta)$ corresponds to the hook
$(\alpha+1,1^\beta)$, while $(d-1,d-2,\dots,1,0|d-1,d-2,\dots,1,0)$ is the
partition $(d^d)$.
\begin{definition}
Let $\CP_\ell$ be the set of partitions of $2\ell$ of the form
$(\alpha_1+1,\dots,\alpha_d+1|\alpha_1,\dots,\alpha_d)$; thus
$\alpha_1+\dots+\alpha_d+d = \ell$ and $\alpha_1>\dots>\alpha_d\ge0$.
\end{definition}

The following plethysm is Ex.\ I.5.10 of Macdonald \cite{Macdonald}:
\begin{equation} \label{plethysm}
\Schur^{(1^\ell)} \circ \Schur^{(2)} = \sum_{\lambda\in\CP_\ell}
\Schur^\lambda .
\end{equation}

\begin{theorem} \label{limit}
The cohomology group $H_n(\L_\HH)(n+\ell)$ is zero except in the following
cases:
\begin{enumerate}
\item $\ell=0$ and $n\ge0$, in which case $H_n(\L_\HH)(n) \cong \HH_n \o
\Schur^{(1^n)}$;
\item $\ell>0$ and $n\ge\ell+2$.
\end{enumerate}
If $\ell>0$ and $n\ge2\ell+2$, we have
$\displaystyle
H_n(\L_\HH)(n+\ell) \cong
\sum_{\substack{\lambda\in\CP_\ell\\n\ge\ell+\alpha_1+1}} \HH_{n-\ell} \o
\Schur^{(1^{n-\ell})+\lambda}$.
\end{theorem}
\begin{proof}
The Chevalley-Eilenberg complex of $\L_\HH(V)$ is bigraded, $\K_{k,\ell} =
\Wedge^k(\HH\o V)\o\Wedge^\ell\bigl(\Schur^2(V)\bigr)$, and since the
differential $\p$ is homogeneous of bidegree $(-2,1)$, the homology is also
bigraded. In terms of this bigrading, we wish to calculate
$H_{n-\ell,\ell}(\L_\HH)$; evidently, this vanishes unless $n\ge\ell$.

The plethysm \eqref{plethysm} implies that
$$
\K_{k,\ell}(\L_\HH)(n) = \sum_{j=0}^{[\frac{k}{2}]}
\sum_{\lambda\in\CP_\ell} \HH_{k-2j} \o \Schur^{(2^j1^{k-2j})} \o
\Schur^\lambda .
$$
We will derive a lower bound for the Laplacian $\Delta$ on each summand.

Given a partition $\lambda\in\CP_\ell$, we calculate that $c_\lambda =
2\ell r + 2\sum_{i=1}^d (\alpha_i + 1) = 2(r+1)\ell$ and
$$
(2^j1^{k-j},\lambda) \le \sum_{i=1}^j (\alpha_i+i+1) + 2\ell \le 3\ell +
\tbinom{j+1}{2} .
$$
On the summand $\HH_{k-2j} \o \Schur^{(2^j1^{k-2j})} \o \Schur^\lambda$, we
have $(r+3)\DD=(r+3)(k+2\ell)$,
\begin{gather*}
\half \Delta_{\GL(V)} \le \half c_{(2^j,1^{k-2j})} + \half c_\lambda +
(2^j1^{k-2j},\lambda) \le \half c_{(2^j,1^{k-2j})} + (r+3)\ell + \ell +
\tbinom{j+1}{2} , \quad \text{and} \\
\begin{aligned}
\Delta_{\SP(\HH)} + c_{(2^j1^{k-2j})} &= \bigl\{ (k-2j)^2 + 2(k-2j) \bigr\}
+ \bigl\{ (k-j)(r-(k-j)+1) + j(r-j+3) \bigr\} \\ &= \half (r+3)k - j(k-j+1)
.
\end{aligned}
\end{gather*}
Combining all of these ingredients, we see that $\Delta \ge
j\bigl(k-\tfrac{3}{2}j+\tfrac{1}{2}\bigr) - \ell$. If $j>0$, the right-hand
side is bounded below by $k-\ell-1$; unless $k\ge2$ and $k\le\ell+1$, our
summand does not contribute to $H_n(\L_\HH)(n+\ell)$. Equivalently,
$n=k+\ell$ must lie in the interval $[\ell+2,2\ell+2]$.

It remains to consider the summands of $\K_{k,\ell}$ with $j=0$; these have
the form
$$
\HH_k \o \sum_{\lambda\in\CP_\ell} \Schur^{(1^k)} \o \Schur^\lambda .
$$
On the summand $\HH_k\o\Schur^{(1^k)+\lambda}$ of
$\HH_k\o\Schur^{(1^k)}\o\Schur^\lambda$, the operator
$\Delta_{\SP(\HH)}+\Delta_{\GL(V)}$ equals
\begin{align*}
k(k+2) + c_{(1^k)} + c_\lambda + 2(1^k,\lambda) &= k(k+2) + k(r-k+1) +
2\ell(r+1) + 2\sum_{i=1}^k \lambda_i \\ &= (k+2\ell)(r+3) -
\sum_{i=k+1}^{\alpha_1+1} \lambda_i ,
\end{align*}
while on all other irreducible components of
$\HH_k\o\Schur^{(1^k)}\o\Schur^\lambda$, it is strictly less. It follows
that the Laplacian can only vanish on the summand
$\HH_k\o\Schur^{(1^k)+\lambda}$, and only at that when $k\ge\alpha_1+1$.
\end{proof}

The following formula illustrates the behaviour of $H_n(\L_\HH)(n+\ell)$
when $n\in[\ell+2,2\ell+1]$
\begin{proposition}
$$
H_n(\L_\HH)(n+1) \cong \Bigl( \HH_{n-1} \o \Schur^{(3,1^{n-2})} \Bigr)
\oplus \begin{cases} \HH_0 \o \Schur^{(4)} , & n=3 , \\ 0 , & n\ne3 .
\end{cases}
$$
\end{proposition}
\begin{proof}
Pieri's formula shows that
\begin{align*}
\K_{n-1,1} &\cong \sum_{j=1}^{[\frac{n+1}{2}]} \HH_{n-2j+1} \o
\Schur^{(2^{j-1},1^{n-2j+1})} \o \Schur^{(2)} \\ &\cong
\sum_{j=1}^{[\frac{n+1}{2}]} \HH_{n-2j+1} \o \Schur^{(2^j,1^{n-2j+1})}
\oplus \sum_{j=1}^{[\frac{n}{2}]} \HH_{n-2j+1} \o
\Schur^{(3,2^{j-1},1^{n-2j})} \\ & \quad \oplus
\sum_{j=2}^{[\frac{n+1}{2}]} \HH_{n-2j+1} \o
\Schur^{(3,2^{j-2},1^{n-2j+2})} \oplus \sum_{j=2}^{[\frac{n+1}{2}]}
\HH_{n-2j+1} \o \Schur^{(4,2^{j-2},1^{n-2j+1})} .
\end{align*}
On these four summands, the operator $\Delta$ equals $j(n-j+2)$,
$j(n-j+3)-n-2$, $j(n-j+1)$ and $j(n-j+2)-n-3$, respectively. Thus, the only
summands on which $\Delta$ vanishes are $\HH_{n-1}\o\Schur^{(3,1^{n-2})}$,
and $\HH_0\o\Schur^{(4)}$.
\end{proof}

The same method may be used in the case $\ell=2$: we obtain
$$
H_n(\L_\HH)(n+2) \cong \Bigl( \HH_{n-2} \o \Schur^{(4,2,1^{n-4})} \Bigr)
\oplus \begin{cases} \HH_1 \o \Schur^{(5,2)} , & n=5 , \\ 0 , & n\ne5
. \end{cases}
$$
Our search for a formula for $H_n(\L_\HH)(n+\ell)$ for all $\ell$ has been
fruitless; nevertheless, it might be of interest to find one.

\end{document}